\documentstyle[12pt]{article} 

\begin{document}

\begin{center}
\bf
\LARGE
Twenty Open Problems \\
in Enumeration of Matchings \\
$ $ \\
\large
\rm
James Propp \\
Mathematical Sciences Research Institute \\
December 9, 1996 \\
(last revised: September 2, 1998)
\end{center}
\normalsize

\bigskip

This document is an exposition of an assortment of open problems
arising from the exact enumeration of (perfect) matchings of finite graphs.
Ten years ago, there were few known results of this kind,
and exact enumeration of matchings
would not have been recognized
as a research topic in its own right.
That situation began to change
in the late '80s with the advent of the Aztec diamond
(in work of Noam Elkies, Greg Kuperberg, Michael Larsen and myself)
and with subsequent study of related problems
by Bo-Yin Yang, William Jockusch, Mihai Ciucu and many others.
Connections with the theory of plane partitions
were brought to the fore by Greg Kuperberg,
who used the matchings viewpoint to prove
a previously conjectural enumeration
of a symmetry class of plane partitions.
Now there are many results in the field,
and even more open problems,
many of which are accessible
to the general combinatorialist.
This article discusses twenty such problems,
some of which are important in their own right,
and others of which are probably merely special cases
of other more general (and more interesting) results.
In most cases, the problem is to find or prove a general formula;
in one case, it is to show that the number of matchings
must be of a certain form (namely, a perfect square),
and in another case, the challenge is to find a combinatorial proof
of a fact whose only known proof is algebraic.

For updates on the status of these problems, see

\begin{center}
{\tt http://www-math.mit.edu/$\sim$propp/progress.ps.gz}.
\end{center}

We begin with problems related to lozenge tilings of hexagons.
A {\bf lozenge\/} is a rhombus of side-length 1
whose internal angles measure 60 and 120 degrees;
all the hexagons we will consider will tacitly have
integer side-lengths and internal angles of 120 degrees.
Every such hexagon $H$ can be dissected into unit equilateral triangles
in a unique way,
and one can use this dissection to define a graph $G$
whose vertices correspond to the triangles
and whose edges correspond to triangles that share an edge;
this is the ``finite honeycomb graph''
dual to the dissection.
It is easy to see
that the tilings of $H$ by lozenges are in one-to-one correspondence
with the perfect matchings of $G$.

The $a,b,c$ semiregular hexagon
is the hexagon whose side lengths are $a,b,c,a,b,c$ respectively.
Lozenge-tilings of this region
are in correspondence with
plane partitions with
at most $a$ rows,
at most $b$ columns,
and no part exceeding $c$.
We may represent such hexagons by means of diagrams like
\begin{verbatim}
                        AVAVAVAVA
                       AVAVAVAVAVA
                      AVAVAVAVAVAVA
                     AVAVAVAVAVAVAVA
                     VAVAVAVAVAVAVAV
                      VAVAVAVAVAVAV
                       VAVAVAVAVAV
                        VAVAVAVAV
\end{verbatim}
where A's and V's represent upward-pointing and downward-pointing triangles,
respectively.

MacMahon showed that the number of such plane partitions is
$$\prod_{i=0}^{a-1} \prod_{j=0}^{b-1} \prod_{k=0}^{c-1}
\frac{i+j+k+2}{i+j+k+1}$$.
(This form of MacMahon's formula is due to Macdonald;
for a short, self-contained proof see section 2 of
``The Shape of a Typical Boxed Plane Partition''
by Henry Cohn, Michael Larsen, and myself,
available at {\tt http://math.wisc.edu/$\sim$propp/shape.ps.gz}).

{\bf Problem 1:} Show that in the $2n-1,2n,2n-1$ semiregular hexagon,
the central location (consisting of the two innermost triangles) is
covered by a lozenge in exactly one-third of the tilings.

(Equivalently: Show that if one chooses a random perfect matching of
the dual graph, the probability that the central edge is contained
in the matching is exactly 1/3.)

\bigskip

The hexagon of side-lengths $n,n+1,n,n+1,n,n+1$ cannot be tiled by
lozenges at all, for in the dissection into unit triangles, the number of 
upward-pointing triangles differs from the number of downward-pointing 
triangles.  However, if one removes the central triangle, one gets a
region that can be tiled, and the sort of numbers one gets for small
values of $n$ are striking.  Here they are, in factored form:

\begin{verbatim}

                             (2)

                                  3
                           (2) (3)

                           5    3
                        (2)  (3)  (5)

                             5    7
                          (2)  (5)

                           2    7    5
                        (2)  (5)  (7)

                        8    3        11
                     (2)  (3)  (5) (7)

                       13    9    11
                    (2)   (3)  (7)   (11)

                       13    18    5     7
                    (2)   (3)   (7)  (11)

                      8    18     13     5
                   (2)  (3)   (11)   (13)

                      2    9     19     11
                   (2)  (3)  (11)   (13)

                   10    3     19     17
                (2)   (3)  (11)   (13)   (17)

                     16     13     23     7
                  (2)   (11)   (13)   (17)

\end{verbatim}
These are similar to the numbers one gets from counting lozenge-tilings
of an $n,n,n,n,n,n$ hexagon, in that the largest prime factor seems to
be bounded by a linear function of $n$, and it ought to be possible to
come up with a conjectural exact formula.  What might be harder is
proving it. 

{\bf Problem 2:} Enumerate the lozenge-tilings of the region obtained
from the $n,n+1,n,n+1,n,n+1$ hexagon by removing the central triangle.

\bigskip

One can also look take a $2n,2n+3,2n,2n+3,2n,2n+3$ hexagon
and make it lozenge-tilable by removing a triangle from the
middle of each of its three long sides, as shown below:

\begin{verbatim}
                       AVAVAVAVAVAVAVAVA
                      AVAVAVAVAVAVAVAVAVA
                     AVAVAVAVAVAVAVAVAVAVA
                    AVAVAVAVAVAVAVAVAVAVAVA
                   AVAVAVAVAVAVAVAVAVAVAVAVA
                   VAVAVAVAVAVAVAVAVAVAVAVAV
                 AVAVAVAVAVAVAVAVAVAVAVAVAVAVA
                AVAVAVAVAVAVAVAVAVAVAVAVAVAVAVA
               AVAVAVAVAVAVAVAVAVAVAVAVAVAVAVAVA
              AVAVAVAVAVAVAVAVAVAVAVAVAVAVAVAVAVA
             AVAVAVAVAVAVAVAVAVAVAVAVAVAVAVAVAVAVA
             VAVAVAVAVAVAVAVAVAVAVAVAVAVAVAVAVAVAV
              VAVAVAVAVAVAVAVAVAVAVAVAVAVAVAVAVAV
               VAVAVAVAVAVAVAVAVAVAVAVAVAVAVAVAV
                VAVAVAVAVAVAVAVAVAVAVAVAVAVAVAV
                 VAVAVAVAVAVAVAVAVAVAVAVAVAVAV
                  VAVAVAVAVAVAVAVAVAVAVAVAVAV
                   VAVAVAVAVAVAVAVAVAVAVAVAV
                    VAVAVAVAVAV VAVAVAVAVAV
\end{verbatim}

Here one obtains an equally tantalizing sequence of factorizations:


\begin{verbatim}

                                1

                              7    2
                           (2)  (7)

                        2    4     4     2
                     (2)  (7)  (11)  (13)

                  10    3    8     2     4     2
               (2)   (3)  (5)  (13)  (17)  (19)

             2    2    2     3     4     4     8     4
          (2)  (5)  (7)  (11)  (13)  (17)  (19)  (23)
\end{verbatim}

{\bf Problem 3:} Enumerate the lozenge-tilings of the region obtained
from the $2n,2n+3,2n,2n+3,2n,2n+3$ hexagon by removing a triangle from
the middle of each of its long sides.

\bigskip

Let us now return to ordinary $a,b,c$ semiregular hexagons.
When $a=b=c$ ($=n$, say), there are not two but six central triangles.
There are two geometrically distinct ways in which we can choose to
remove an upward-pointing triangle and downward-pointing triangle
from these six, according to whether the triangles are opposite or adjacent:
\begin{verbatim}
                 AVAVAVA                 AVAVAVA
                AVAVAVAVA               AVAVAVAVA
               AVAVA AVAVA             AVAV VAVAVA
               VAVAV VAVAV             VAVA AVAVAV
                VAVAVAVAV               VAVAVAVAV
                 VAVAVAV                 VAVAVAV
\end{verbatim}
Such regions may be called ``holey hexagons'' of two different kinds.
In the former case, the number of tilings of the holey hexagon is a 
nice round number (in the sense that, like the numbers tabulated above for 
Problems 2 and 3, it has small prime factors).  In the latter case, the 
number of tilings is not round.  Note, however, that in the latter case,
the number of tilings of the holey hexagon divided by the number of tilings 
of the unaltered hexagon (given to us by MacMahon's formula) is equal to 
the probability that a random lozenge tiling of the hexagon contains a 
lozenge that covers these two triangles, and tends to 1/3 for large $n$.
Following this clue, we examine the difference between the aforementioned 
probability (with its messy, un-round numerator) and the number 1/3.  The 
result is a fraction in which the numerator is now a nice round number.  
So, in both cases, we may have reason to think that there is an exact formula.

{\bf Problem 4:} Determine the number of lozenge-tilings of a
regular hexagon from which two of its innermost unit triangles
(one upward-pointing and one downward-pointing) have been removed.

\bigskip

At this point, I must digress and explain how I did the exploratory work
that indicated that all these numbers are ``nice''.  Greg Kuperberg
wrote a program called {\tt dommaple}, which (with enhancements by 
David Wilson and myself) became a program called {\tt vaxmaple}.  One
can feed {\tt vaxmaple} an ASCII file of V's and A's like the ones
shown above, and it will output Maple code which, if piped to {\tt maple},
will output the number of tilings of the region.  Moreover, as we will
see below, {\tt vaxmaple} can also count domino-tilings of regions
(indeed, {\tt dommaple} can do this too; the main difference between the
two programs is that {\tt vaxmaple} can handle lozenges as well).
Then there is {\tt vaxmacs}, which provides a way to do all this 
interactively in real-time; interested readers with access to 
the World Wide Web can obtain copies of both programs via
{\tt http://math.wisc.edu/~propp/software.html}.

The main point I want to make here is that these programs take advantage
of a result of Percus (based on work of Kasteleyn) that says that the 
number of matchings of a bipartite planar graph can always be calculated 
as the absolute value of the determinant of a modified adjacency matrix $K$
consisting of $0$'s, $1$'s, and
$-1$'s, in which the rows correspond to vertices in one component of
the bipartition and the columns correspond to vertices in the other component.
In the case of lozenge tilings of hexagons and the associated matchings,
it turns out that there is no need to modify signs of entries; the ordinary
adjacency matrix will do.

Unpublished work of Greg Kuperberg shows that 
when row-reduction and column-reduction are
systematically applied to the Kasteleyn matrix of 
an $a,b,c$ semiregular hexagon, 
one can obtain the $b$-by-$b$ Carlitz matrix whose $i,j$th entry is 
$a+c \choose a+i-j$.
(This matrix can also be recognized as the Gessel-Viennot matrix
that arises from interpreting each tiling
as a family of non-intersecting lattice paths;
Horst Sachs and his colleagues noticed this
a number of years ago.)
Such reductions do not affect the determinant,
so we have a pleasing way of understanding the relationship
between the Kasteleyn-Percus matrix method
and the Gessel-Viennot lattice-path method.
However, one can also verify that the reductions
do not affect the {\it cokernel} of the matrix, either.
On the other hand, the cokernel of the Kasteleyn matrix
for the $a,b,c$ hexagon is clearly invariant under
permuting $a$, $b$, and $c$.
This gives rise to three different Carlitz matrices
that non-trivially have the same cokernel.
E.g., if $c=1$,
then one gets an $a$-by-$a$ matrix and a $b$-by-$b$ matrix
that both have the same cokernel,
which can be computed by noticing that the third Carlitz matrix of the trio
is just a 1-by-1 matrix
whose sole entry is (plus or minus) a binomial coefficient.
In this special case,
the cokernel is just a cyclic group.

Greg Kuperberg poses the challenge:

{\bf Problem 5:} Determine the cokernel of the Carlitz matrices,
or equivalently of the Kasteleyn matrices of $a,b,c$ hexagons, and
if possible find a way to interpret the cokernel in terms of
the tilings.

He points out that in the case $a=b=c=2$, one gets a non-cyclic group
as the cokernel.

\bigskip

Digressing for a moment from the topic of lozenge tilings,
I point out that
in general, the Kasteleyn matrix $K$ is not canonically defined,
in the sense that there may be many ways of modifying the
signs of certain entries of the bipartite adjacency matrix
of a graph so that all non-zero contributions to the
determinant have the same sign.  Thus, one should not expect
the eigenvalues of $K$ to possess combinatorial significance.
However, the spectrum of $K$ times its adjoint $K^*$ is independent
of which Kasteleyn matrix $K$ one chooses (independently shown
by David Wilson and Horst Sachs).  Thus, it is natural to ask:

{\bf Problem 6:} What is the significance of the spectrum of
$K K^*$, where $K$ is any Kasteleyn matrix associated with a
bipartite planar graph?

\bigskip

Returning now to lozenge tilings, or equivalently, perfect matchings
of finite subgraphs of the infinite honeycomb,
consider the hexagon graph with $a=b=c=2$:

\pagebreak

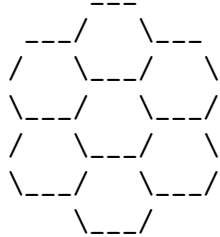
\begin{verbatim}
                           ___
                       ___/   \___
                      /   \___/   \
                      \___/   \___/
                      /   \___/   \
                      \___/   \___/
                          \___/   

\end{verbatim}
This is the graph whose 20 perfect matchings correspond to the 20 tilings
of the regular hexagon of side 2 by rhombuses of side 1.  If we just look
at the probability of each individual horizontal edge belonging to a
matching chosen uniformly at random (``edge-probabilities''), we get
\begin{verbatim}
                            .7
                        .3      .3
                            .3     
                        .4      .4 
                            .3     
                        .3      .3 
                            .7    
\end{verbatim}
Now let us look at this table of numbers as if it described a distribution
of mass.  If we assign the three columns $x$-coordinates $-1$ through 1, 
we find that the weighted sum of the squares of the $x$-coordinates is equal to
$(.3+.4+.3) (-1)^2  +  (.7+.3+.7+.3) (0)^2  +  (.3+.4+.3) (1)^2$, or
$(1.0) (-1)^2 + (2.0) (0)^2 + (1.0) (1)^2$, or 2.  If we assign the seven 
rows $y$-coordinates $-3$ through 3, we find that the weighted sum of the 
squares of the $y$-coordinates is equal to
$(.7) (-3)^2  +  (.6) (-2)^2  +  (.3) (-1)^2  +  (.8) (0)^2
+  (.3) ( 1)^2  +  (.6) ( 2)^2  +  (.7) (3)^2$, or 20.
You can do a similar (but even easier) calculation yourself for the case
$a=b=c=1$, to see that the ``moments of inertia'' of the horizontal
edge-probabilities around the vertical and horizontal axes are 0 and
1, respectively.  Using {\tt vaxmaple} to study the case $a=b=c=n$ for
larger values of $n$, I find that the moment of inertia about the vertical 
axis goes like
$$0, 2, 12, 40, 100, ...$$
and the moment of inertia about the horizontal axis goes like
$$1, 20, 93, 296, 725, ... .$$
It is easy to show that the former numbers 
are given in general by the polynomial
$(n^4-n^2)/6$.
The latter numbers are subtler;
they are not given by a polynomial of degree 4,
though it is noteworthy that the $n$th term
is an integer divisible by $n$,
at least for the first few values of $n$.

{\bf Problem 7:} Find the ``vertical moments of inertia'' for the mass on
edges arising from edge-probabilities for random matchings of the 
$a,b,c$ honeycomb graph.

\bigskip

Now let us turn from lozenge-tiling problems to domino-tiling problems.
A {\bf domino} is a 1-by-2 or 2-by-1 rectangle.
Although lozenge-tilings (in the guise of constrained plane partitions)
were studied first, it was really the study of domino tilings in Aztec
diamonds that gave current work on enumeration of matchings its current 
impetus.  Here is the Aztec diamond of order 5:
\begin{verbatim}
                              XX
                             XXXX
                            XXXXXX
                           XXXXXXXX
                          XXXXXXXXXX
                          XXXXXXXXXX
                           XXXXXXXX
                            XXXXXX
                             XXXX
                              XX
\end{verbatim}
(An X represents a 1-by-1 square.)
A tiling of such a region by dominos is equivalent to a perfect matching
of a certain (dual) subgraph of the infinite square grid.
This grid is bipartite, and it is convenient to color its vertices
alternately black and white;
equivalently, it is convenient to color the 1-by-1 squares alternately 
black and white, so that every domino contains 
one 1-by-1 square of each color.
Elkies, Kuperberg, Larsen, and Propp showed that the number of domino-tilings
of such a region is $2^{n(n+1)/2}$ (where $2n$ is the number of rows),
and Ionescu later proved an exact formula
(originally conjectured by Jockusch) for the number of tilings of regions like

\pagebreak

\begin{verbatim}
                              XX
                             XXXX
                            XXXXXX
                           XXXXXXXX
                          XXXX XXXXX
                          XXXX XXXXX
                           XXXXXXXX
                            XXXXXX
                             XXXX
                              XX
\end{verbatim}
in which two innermost squares of opposite color have been removed.

Now suppose one removes two squares from the middle of an Aztec diamond of 
order $n$ in the following way:
\begin{verbatim}
                              XX
                             XXXX
                            XXXXXX
                           XXXX XXX
                          XXXXXXXXXX
                          XXXX XXXXX
                           XXXXXXXX
                            XXXXXX
                             XXXX
                              XX
\end{verbatim}
(The two squares removed are a knight's-move apart, and subject to that
constraint, they are as close to being in the middle as they can be.
Up to symmetries of the square, there is only one way of doing this.)
Then numbers of tilings one gets are as follows (for $n = 2$ through 10):

\pagebreak

\begin{verbatim}

                           (2)

                              3
                           (2)

                            5
                         (2)  (5)

                           9    2
                        (2)  (3)

                           17
                        (2)   (3)

                           22    2
                        (2)   (3)

                        24    2
                     (2)   (3)  (73)

                      31    2    2
                   (2)   (3)  (5)  (11)

                         47    2
                      (2)   (3)  (5)
\end{verbatim}
Note that only the presence of the large prime factor 73 makes one doubt
that there is a general formula; the other prime factors are reassuringly
small.  Further data might make it clear what that 73 is doing there.

{\bf Problem 8:} Count the domino tilings of an Aztec diamond from which
two close-to-central squares, related by a knight's move, have been
deleted.

\bigskip

One can also look at ``Aztec rectangles'' from which squares have been
removed so as to restore the balance between black and white squares
(a necessary condition for tilability).  For instance, one can remove
the central square from an $a$-by-$b$ Aztec rectangle in which $a$
and $b$ differ by 1, with the larger of $a,b$ odd:
\begin{verbatim}
                              XX
                             XXXX
                            XXXXXX
                           XXXXXXXX
                           XXXX XXXX
                            XXXXXXXX
                             XXXXXX
                              XXXX
                               XX
\end{verbatim}

{\bf Problem 9:} Find a formula for the number of domino tilings of a
$2n$-by-$(2n+1)$ Aztec rectangle with its central square removed.

\bigskip

What about $(2n-1)$-by-$2n$ rectangles?  For these regions, removing
the central square does not make the region tilable.  However, if
one removes any one of the four squares adjacent to the middle square,
one obtains a region that is tilable, and moreover, for this region
the number of tilings appears to be a nice round number. 

{\bf Problem 10:} Find a formula for the number of domino tilings of a
$(2n-1)$-by-$2n$ Aztec rectangle with a square adjoining the central
square removed.

\bigskip

At this point, readers who are unfamiliar with the literature may be
wondering why $m$-by-$n$ rectangles haven't come into the story.
Indeed, one of the surprising facts of life in the study of enumeration
of matchings is that Aztec diamonds and their kin have (so far) been
much more fertile ground for exact combinatorics that the seemingly
more natural rectangles.  There are, however,
a few cases I know of in which something rather nice turns up.  One
is the problem of Ira Gessel that appears as Problem 20 in this document.
Another is the work done by Jockusch and, later, Ciucu on why the number 
of domino tilings of the square is always either a perfect square or
twice a perfect square.  In the spirit of the work of Jockusch and
Ciucu, I offer here a problem based on Lior Pachter's observation
that the region

\pagebreak

\begin{verbatim}
                     XXXXXXXXXXXXXXXX
                     XXXXXXXXXXXXXXXX
                     XXXXXXXXXXXXXXXX
                     XXXXXXXXXXXXXXXX
                     XXXXXXXXXXXXXXXX
                     XXXXXXXXXXXXXXXX
                     XXXXXXXXXXXXXXXX
                     XXXXXXXXXXXXXXXX  
                     XXXXXXX  XXXXXXX
                     XXXXXX  XXXXXXXX
                     XXXXX  XXXXXXXXX
                     XXXX  XXXXXXXXXX
                     XXX  XXXXXXXXXXX
                     XX  XXXXXXXXXXXX
                     X  XXXXXXXXXXXXX
                       XXXXXXXXXXXXXX
\end{verbatim}
(8 dominos removed from a 16x16 square) has exactly one tiling.  
What if we make the intrusion half as long, as in the following picture?
\begin{verbatim}
                     XXXXXXXXXXXXXXXX
                     XXXXXXXXXXXXXXXX
                     XXXXXXXXXXXXXXXX
                     XXXXXXXXXXXXXXXX
                     XXXXXXXXXXXXXXXX
                     XXXXXXXXXXXXXXXX
                     XXXXXXXXXXXXXXXX
                     XXXXXXXXXXXXXXXX
                     XXXXXXXXXXXXXXXX
                     XXXXXXXXXXXXXXXX
                     XXXXXXXXXXXXXXXX
                     XXXXXXXXXXXXXXXX
                     XXX  XXXXXXXXXXX
                     XX  XXXXXXXXXXXX
                     X  XXXXXXXXXXXXX
                       XXXXXXXXXXXXXX
\end{verbatim}
That is, we take a $2n$-by-$2n$ square (with $n$ even) and remove $n/2$ 
dominos from it, in a partial zig-zag pattern that starts from the corner.
Here are the numbers we get, in factored form, for $n=2,4,6,8,10$:
\begin{verbatim}

                                   2
                            (2) (3)

                           2    6     2
                        (2)  (3)  (13)

                     3    2    4    2       2
                  (2)  (3)  (5)  (7)  (3187)

                      4           2        2
                   (2)  (11771899)  (27487)

                     5                      2
                  (2)  (2534588575976069659)
\end{verbatim}
The factors are ugly, but the exponents are nice: we get $2^{n/2}$ times 
an odd square.

Perhaps this is a special case of a two-parameter fact that says that
you can take an intrusion of length $m$ in a $2n$-by-$2n$ square and the
number of tilings of the resulting region will always be a square or
twice a square.

{\bf Problem 11:} What is going on with ``intruded Aztec diamonds''?
In particular, why is the number of tilings so square-ish?

\bigskip

Let's now get back to those Kasteleyn matrices we discussed earlier.
Work of Rick Kenyon and David Wilson has shown that the {\it inverses}
of these matrices are loaded with combinatorial information, so it
would be nice to get our hands on them.  Unfortunately, there are
a lot of non-zero entries in the inverse-matrices.  (Recall that
the Kasteleyn matrices themselves, being nothing more than adjacency
matrices in which some of the 1's have been strategically replaced
by $-1$'s, are sparse; their inverses, however, tend to have most
if not all of their entries non-zero.)  Nonetheless, some exploratory
numerology leaves room for hope that this is do-able.

Consider the Kasteleyn matrix $K_n$ for the Aztec diamond of order $n$,
in which every other vertical domino has its sign flipped (that is,
the corresponding 1's in the bipartite adjacency matrix are
replaced by $-1$'s). 

{\bf Problem 12:} Show that the sum of the entries of the matrix inverse
of $K_n$ is $(n-1)(n+3)/2 - 2^{n-1} + 2$.

(This formula works for $n=1$ through $n=8$.)

I should mention in this connection that Greg Kuperberg has some
high-tech ruminations on the inverses of Kasteleyn matrices, and
there is a chance that representation-theory methods will be useful
here.

\bigskip

Now we turn to a class of regions I call ``pillows''
on account of their agreeably lumpy shape.
Here is a ``0 mod 4'' pillow of ``order 5'':
\begin{verbatim}

                            XXXX
                         XXXXXXXX
                      XXXXXXXXXXXX
                   XXXXXXXXXXXXXXXX
                XXXXXXXXXXXXXXXXXXXX
                XXXXXXXXXXXXXXXXXXXX
                 XXXXXXXXXXXXXXXX
                  XXXXXXXXXXXX
                   XXXXXXXX
                    XXXX
\end{verbatim}
And here is a ``2 mod 4'' pillow of ``order 7'':

\pagebreak

\begin{verbatim}
                              XX
                           XXXXXX
                        XXXXXXXXXX
                     XXXXXXXXXXXXXX
                  XXXXXXXXXXXXXXXXXX
               XXXXXXXXXXXXXXXXXXXXXX
            XXXXXXXXXXXXXXXXXXXXXXXXXX
            XXXXXXXXXXXXXXXXXXXXXXXXXX
             XXXXXXXXXXXXXXXXXXXXXX
              XXXXXXXXXXXXXXXXXX
               XXXXXXXXXXXXXX
                XXXXXXXXXX
                 XXXXXX
                  XX
\end{verbatim}
It turns out (empirically) that the number of 0-mod-4 pillows of order $n$
is a perfect square times the coefficient of $x^n$ in the Taylor expansion
of $(5+3x+x^2-x^3)/(1-2x-2x^2-2x^3+x^4)$.  Similarly, it appears that the
number of 2-mod-4 pillows of order $n$ is a perfect square times the
coefficient of $x^n$ in the Taylor expansion of
$(5+6x+3x^2-2x^3)/(1-2x-2x^2-2x^3+x^4)$.  (If you're wondering about
``odd pillows'', I should mention that there's a nice formula for the
number of tilings, but it isn't an interesting result, because an odd pillow
splits up into many small non-communicating sub-regions such that a tiling 
of the whole region corresponds to a choice of tiling on each of the 
sub-regions.)

{\bf Problem 13:} Find a general formula for the number of domino tilings
of even pillows.

\bigskip

Jockusch looked at the order-$n$ Aztec diamond with a 2-by-2
hole in the center, for small values of $n$;
he came up with a conjecture for the number of domino tilings,
subsequently proved by Ionescu.
One way to generalize this is to make the hole larger,
as was suggested by Doug Zare.
Here's what David Wilson reported on this subject
on October 15, 1996 (via e-mail):

\begin{verbatim}
We all know the formula for the number of tilings that the 
Aztec diamond has.  Doug Zare asked how many tilings there 
are in an Aztec diamond with an Aztec diamond deleted from 
it.  Let us define the Aztec window with outer order y and 
inner order x to be the Aztec diamond of order y with an 
Aztec diamond of order x deleted from its center.  For 
example, this is the Aztec window with orders 8 and 2:
\end{verbatim}


\begin{verbatim}

                            XX       
                           XXXX      
                          XXXXXX     
                         XXXXXXXX    
                        XXXXXXXXXX   
                       XXXXXXXXXXXX  
                      XXXXXX  XXXXXX 
                     XXXXXX    XXXXXX
                     XXXXXX    XXXXXX
                      XXXXXX  XXXXXX 
                       XXXXXXXXXXXX  
                        XXXXXXXXXX   
                         XXXXXXXX    
                          XXXXXX     
                           XXXX      
                            XX       


There are a number of interesting patterns that show up 
when we count tilings of Aztec windows.  For one thing, 
if w is a fixed even number, and y = x+w, then for any w 
the number of tilings appears to be a polynomial in x.  
(When w is odd, and x is large enough, there are no tilings.)  
For w=6, the polynomial is
           8          7           6            5            4
     8192 x  + 98304 x  + 573440 x  + 2064384 x  + 4988928 x  

                     3          2
          + 8257536 x  9175040 x  + 6291456 x + 2097152.

Substituting x=2, the above region has 314703872 tilings.
This isn't just some random polynomial.  It can be rewritten
as

        (2)^17*x^4   1/2*x+7/8   (x+3/2)^2

where these three polynomials get composed.
For the above example, we evaluate (2+3/2)^2 = 49/4,
1/2*49/4 + 7/8 = 7, 2^(17)*7^4 = 314703872 tilings.
For all integer values of x, by the time the second polynomial
is applied, the result is an integer.

Is this decomposition of the 6th Aztec window polynomial a 
fluke?  Of course not!  In general the rightmost polynomial 
is (x+w/4)^2, and the leftmost polynomial is either a perfect 
square, twice a fourth power, or half a fourth power, depending 
on w mod 8.  A pattern for the middle polynomial however 
eludes me.

We all know that the Aztec window polynomials will be squareish
because of symmetry, but why would they be quarticish half the
time, but only perfect squares the other half of the time?
The rightmost polynomial in the decomposition is equivalent to
saying that when the polynomials are expressed in terms of 
(3*x+y)/4 rather than x, there are no odd degree terms.  
I gave these polynomials to a computer-algebra person who said 
that he could find functional decompositions of polynomials, 
but it turns out that no-one implemented his algorithm, so I 
don't know whether or not the polynomials decompose further.

Using the old version of vaxmacs I was able to determine the 
polynomials for values of w up to 14, with the new version 
(which is now installed at MSRI) it was not too hard to compute 
the polynomials for all w up to 34.

Does anybody see a pattern to these polynomials?  Or how to 
prove the above observations?  The first few polynomials are 
given below, and are normalized so that the left polynomial has 
the largest constant factor consistent with the composition of 
the middle and right polynomials being integer-valued for integer 
x.


(2)^3*x^4   1   (x+1/2)^2           [w=2]

(2)^8*x^2   x+1   (x+1)^2           [w=4]

(2)^17*x^4   1/2*x+7/8   (x+3/2)^2  [w=6]

(2)^28*x^2   1/144*x^4+7/72*x^3+41/144*x^2+11/18*x+1   (x+2)^2

(2)^43*x^4   1/144*x^3+61/576*x^2+451/2304*x+967/1024   (x+5/2)^2

\end{verbatim}

(David e-mailed much more data, but I'll omit details here.)

{\bf Problem 14:} Find a general formula for the number of domino tilings
of Aztec windows.

Even an argument explaining why the number of tilings for windows
of inner order $x$ and outer order $x+w$  should be given by a polynomial
in $x$ (for each fixed $w$) would constitute progress!

\bigskip

Now we come to some problems involving tiling that fit neither the
domino-tiling nor the lozenge-tiling paradigm.
Here the more general picture is that we have some periodic dissection 
of the plane by polygons,
such that an even number of polygons meet at each vertex,
allowing us to color the polygons alternately black or white.
We then make a clever choice of a finite region $R$
composed of equal numbers of black and white polygons,
and we look at the number of ``diform'' tilings of the region,
where a {\bf diform} is the union of two polygonal cells
that share an edge.
In the case of domino-tilings, the underlying dissection of the infinite plane
is the tiling by squares, 4 around each vertex;
in the case of lozenge-tilings, the underlying dissection of the infinite plane
is the tiling by equilateral triangles, 6 around each vertex.

Other sorts of periodic dissections have already played a role in
the theory of enumeration of matchings.
For instance, there is a tiling of the plane by isosceles right triangles
associated with a discrete reflection group in the plane;
in this case, the right choice of $R$
gives us a region that can be tiled in
$5^{n^2}$ ways.
Similarly, in the tiling of the plane by triangles
that comes from a 30 degree, 60 degree, 90 degree right triangle
by repeatedly reflecting it in its edges
gives rise to a tiling problem
in which powers of 13 occur.
I cannot include the pictures here,
but I will say that one key feature of these regions $R$
is revealed by looking at the colors of those polygons in the dissection
that share an edge with the border of $R$.
One sees that the border splits up into four long stretches
such that along each stretch,
all the polygons that touch the border have the same color.

One case that has not yet been settled
is the case that arises from a rather symmetric
dissection of the plane into 
equilateral triangles, squares, and regular hexagons,
with 4 polygons meeting at each vertex.
Empirically, one finds that
the number of diform tilings is
$2^{n(n+1)}$.

{\bf Problem 15:} Prove that for the Aztec-type regions in the dissection
of the plane into triangles, squares, and hexagons, the number of tilings
of the region of order $n$ is $2^{n(n+1)}$.

(See {\tt http://www-math.mit.edu/$\sim$propp/dragon.ps}
for a picture of one of these tilings.)

\bigskip

One way to get a new dissection of the plane from an old one is to refine it.
For instance, starting from the dissection of the plane into squares,
one can draw in every $k$th southwest-to-northeast diagonal.  When
$k$ is 1, this is just a distortion of the dissection of the plane into 
equilateral triangles.  When $k$ is 2, this is a dissection that leads to 
finite regions for which the number of diform tilings is a known power of 2
(thanks to a theorem of Chris Douglas).  But what about $k=3$ and higher?

For instance, we have the roughly hexagonal region shown at the top of
the next page
(a union of square and triangular pieces, with certain boundary vertices
marked with a ``.''\ so as to bring out the large-scale 2,3,2,2,3,2 hexagonal 
structure more clearly);
it has 17920 tilings, where 17920 is $2^9 \cdot 5 \cdot 7$.  More generally,
if one takes an $a,b,c$ quasi-hexagon, one finds that one gets a large
power of 2 times a product of powers of odd primes in which all the
primes are fairly small (and their exponents are too).

{\bf Problem 16:} Find a formula for the number of diform tilings in the
$a,b,c$ quasi-hexagon in the dissection of the plane that arises from
slicing the dissection into squares along every third upward-sloping
diagonal.

I should mention that one reason for my special interest in Problem 16
is that it seems to be a genuine hybrid of domino tilings of Aztec diamonds
and lozenge tilings of hexagons.  I should also mention that I think
the problem will yield to some generalization of the method of 
subgraph substitution that has already been of such
great use in enumeration of matchings of graphs.


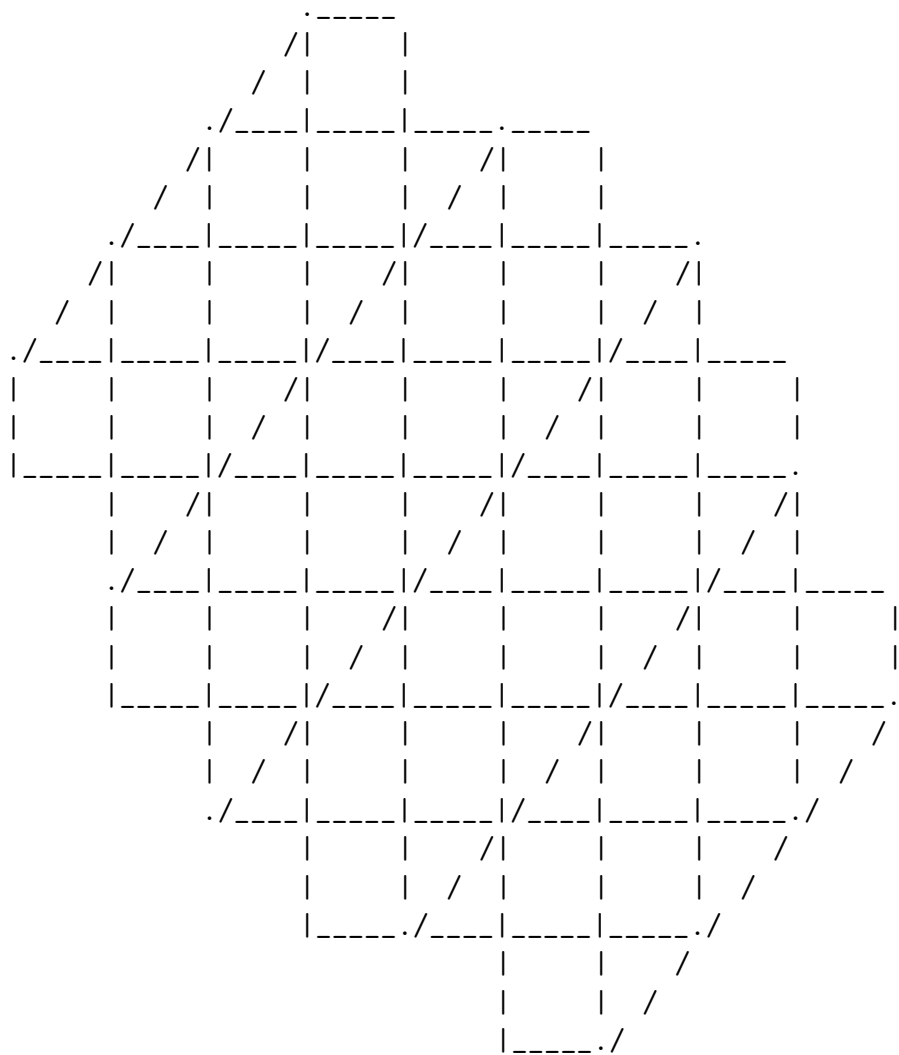
\begin{verbatim}
                        ._____
                       /|     |     
                     /  |     |     
                  ./____|_____|_____._____
                 /|     |     |    /|     |
               /  |     |     |  /  |     |
            ./____|_____|_____|/____|_____|_____.
           /|     |     |    /|     |     |    /|
         /  |     |     |  /  |     |     |  /  |
      ./____|_____|_____|/____|_____|_____|/____|_____ 
      |     |     |    /|     |     |    /|     |     |
      |     |     |  /  |     |     |  /  |     |     |
      |_____|_____|/____|_____|_____|/____|_____|_____.
            |    /|     |     |    /|     |     |    /|
            |  /  |     |     |  /  |     |     |  /  |
            ./____|_____|_____|/____|_____|_____|/____|_____
            |     |     |    /|     |     |    /|     |     |
            |     |     |  /  |     |     |  /  |     |     |
            |_____|_____|/____|_____|_____|/____|_____|_____.
                  |    /|     |     |    /|     |     |    /
                  |  /  |     |     |  /  |     |     |  /
                  ./____|_____|_____|/____|_____|_____./
                        |     |    /|     |     |    / 
                        |     |  /  |     |     |  /   
                        |_____./____|_____|_____./
                                    |     |    / 
                                    |     |  /
                                    |_____./
\end{verbatim}

I will not go into great detail here on this method (the interested
reader can examine
{\tt http://www-math.mit.edu/$\sim$propp/}
{\tt fpsac96.ps.gz}).
Here is an overview:
One studies not graphs but weighted graphs,
with weights assigned to edges,
and one does weighted enumeration of perfect matchings,
where the weight of a matching is the product
of the weights of the constituent edges.
One then looks at local substitutions with a graph that
preserve the sum of the weights of the matchings,
or more generally,
multiply the sum of the weights of the matchings
by some predictable factor.
Then the problem of weight-enumerating matchings of one graph
reduces to the problem of weight-enumerating matchings
of another (hopefully simpler) graph.
Iterating this procedure,
one can often eventually reduce the graph
to something one already understands.
I am confident that this will apply to Problem 16,
yielding a reduction to the standard $a,b,c$ semiregular hexagon.

\bigskip

Problems 15 and 16 are just two instances of a broad class of problems
arising from periodic graphs in the plane.  A unified understanding
of this class of problems has begun to emerge, by way of subgraph
substitution.  The most important open problem connected with this class
of results is the following:

{\bf Problem 17:} Characterize those local substitutions that have 
a predictable effect on the weighted sum of matchings of a graph.

The most useful local substitution so far has been
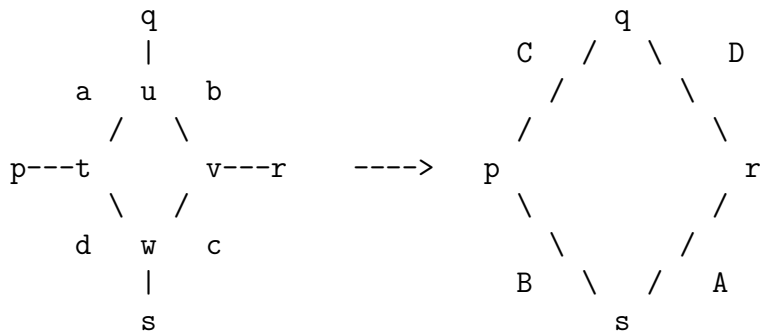
\begin{verbatim}

             q                            q
             |                      C   /   \    D
         a   u   b                    /       \
           /   \                    /           \
     p---t       v---r    ---->   p               r
           \   /                    \           /
         d   w   c                    \       /
             |                      B   \   /   A
             s                            s
\end{verbatim}
(where unmarked edges have weight 1
and where $A,B,C,D$ are obtained from $a,b,c,d$
by dividing by $ac+bd$),
but Rick Kenyon's substitution
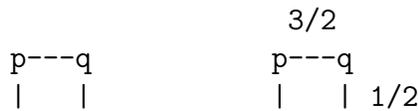
\begin{verbatim}
                                   3/2
                  p---q           p---q
                  |   |           |   | 1/2
                  r---s   ---->   r---s
                  |   |         2 |   | 
                  t---u           t   u
\end{verbatim}
has also been of use
(where vertices $p,q,t,u$ may be attached elsewhere in the graph, 
but not the vertices $r,s$).

\bigskip

Up till now we have been dealing exclusively with bipartite planar graphs.
However, it is possible that there exists rich combinatorics involving
other sorts of graphs.

For instance, one can look at the triangle graph of order $n$:
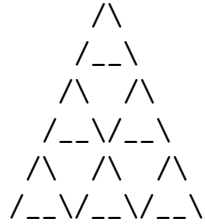
\begin{verbatim}
                            /\
                           /__\
                          /\  /\
                         /__\/__\
                        /\  /\  /\
                       /__\/__\/__\
\end{verbatim}
This particular graph has 6 matchings; we write $M(4) = 6$.
More generally, we let $M(n)$ denote the number of matchings of the
triangular graph whose longest row contains 6 vertices.  
When $n$ is 1 or 2 mod 4, the graph has an odd number of vertices
and $M(n)$ is 0; hence let us only consider the cases in which
$n$ is 0 or 3 mod 4.
Here are the
first few values of $M(n)$, expressed in factored form:
$2$,
$2 \cdot 3$,
$2 \cdot 2 \cdot 3 \cdot 3 \cdot 61$,
$2 \cdot 2 \cdot 11 \cdot 29 \cdot 29$,
$2^3 \cdot 3^3 \cdot 5^2 \cdot 7^2 \cdot 19 \cdot 461$,
$2^3 \cdot 5^2 \cdot 37^2 \cdot 41 \cdot 139^2$,
$2^4 \cdot 73 \cdot 149 \cdot 757 \cdot 33721 \cdot 523657$,
$2^4 \cdot 3^8 \cdot 17 \cdot 37^2 \cdot 703459^2$,
\dots.
It is interesting that $M(n)$ seems to be divisible by
$2^{\lfloor (n+1)/4 \rfloor}$
but no higher power of 2;
it is also interesting that when we divide by this power of 2,
in the case where $n$ is a multiple of 4,
the quotient we get, in addition to being odd,
is a perfect square times a small number
$(3, 11, 41, 17, \dots)$.

{\bf Problem 18:} How many perfect matchings does the triangle graph
of order $n$ have?

\bigskip

One can also look at graphs that are bipartite but not planar.
A natural example is the $n$-cube (that is, the $n$-dimensional cube
with all sides of length 2).  It has been shown that the number
of perfect matchings of the $n$-cube goes like $1$, $2$, $9=3^2$,
$272=16 \cdot 17$,
$589185=3^2 \cdot 5 \cdot 13093$, $\dots$.

{\bf Problem 19:} Find a formula for the number of perfect matchings
of the $n$-cube.

(This may be intractable; after all, the graph has exponentially
many vertices.)

\bigskip

Finally, let's go to a problem involving domino tilings of rectangles,
submitted by Ira Gessel (what follows are his words):



\catcode `!=11

\newdimen\squaresize 
\newdimen\thickness 
\newdimen\Thickness
\newdimen\ll! \newdimen \uu! \newdimen\dd!
\newdimen \rr! \newdimen \temp!

\def\sq!#1#2#3#4#5{%
\ll!=#1 \uu!=#2 \dd!=#3 \rr!=#4
\setbox0=\hbox{%
 \temp!=\squaresize\advance\temp! by .5\uu!
 \rlap{\kern -.5\ll! 
 \vbox{\hrule height \temp! width#1 depth .5\dd!}}%
%
 \temp!=\squaresize\advance\temp! by -.5\uu!  
 \rlap{\raise\temp! 
 \vbox{\hrule height #2 width \squaresize}}%
%
 \rlap{\raise -.5\dd!
 \vbox{\hrule height #3 width \squaresize}}%
%
 \temp!=\squaresize\advance\temp! by .5\uu!
 \rlap{\kern \squaresize \kern-.5\rr! 
 \vbox{\hrule height \temp! width#4 depth .5\dd!}}%
%
 \rlap{\kern .5\squaresize\raise .5\squaresize
 \vbox to 0pt{\vss\hbox to 0pt{\hss $#5$\hss}\vss}}%
}
 \ht0=0pt \dp0=0pt \box0
}

\def\vsq!#1#2#3#4#5\endvsq!{\vbox 
  to \squaresize{\hrule width\squaresize height 0pt%
\vss\sq!{#1}{#2}{#3}{#4}{#5}}}

\newdimen \LL! \newdimen \UU! \newdimen \DD! \newdimen \RR!

\def\vvsq!{\futurelet\next\vvvsq!}
\def\vvvsq!{\relax
  \ifx     \next l\LL!=\Thickness \let\continue!=\skipnexttoken!
  \else\ifx\next u\UU!=\Thickness \let\continue!=\skipnexttoken!
  \else\ifx\next d\DD!=\Thickness \let\continue!=\skipnexttoken!
  \else\ifx\next r\RR!=\Thickness \let\continue!=\skipnexttoken!
  \else\ifx\next P\let\continue!=\place!
  \else\def\continue!{\vsq!\LL!\UU!\DD!\RR!}%
  \fi\fi\fi\fi\fi 
  \continue!}

\def\skipnexttoken!#1{\vvsq!}

\def\place! P#1#2#3{%
\rlap{\kern.5\squaresize\temp!=.5\squaresize\kern#1\temp!
  \temp!=\squaresize 
  \advance\temp! by #2\squaresize \temp!=.5\temp!
  \raise\temp!\vbox 
   to 0pt{\vss\hbox to 0pt{\hss$#3$\hss}\vss}}\vvsq!}

\def\Young#1{\LL!=\thickness \UU!=\thickness
 \DD! = \thickness \RR! = \thickness
\vbox{\smallskip\offinterlineskip
\halign{&\vvsq! ## \endvsq!\cr #1}}}

\def\Youngt#1{\LL!=\thickness \UU!=
  \thickness \DD! = \thickness \RR! = \thickness
\vtop{\offinterlineskip
\halign{&\vvsq! ## \endvsq!\cr #1}}}

\def\blank{\omit\hskip\squaresize}
\catcode `!=12

\thickness=.4pt

We consider dimer coverings of an $m\times n$ rectangle,
with $m$ and $n$ even. We assign a vertical domino from row
$i$ to row $i+1$ the weight
$\sqrt {y_i}$ and a horizontal domino from column $j$ to
column $j+1$ the weight $\sqrt {x_j}$. For example, the
covering
$$\Thickness=0pt
\squaresize=30pt
\Young{dP0{-1}{\sqrt{y_1}}&rP10
  {\sqrt{x_2}}&l&dP0{-1}
  {\sqrt{y_1}}&
  rP10{\sqrt{x_5}}&l&rP10{\sqrt{x_7}}
  &l&dP0{-1}{\sqrt{y_1}}&dP0{-1}{\sqrt{y_1}}\cr
  u&rP10{\sqrt{x_2}}&l&u&rP10
  {\sqrt{x_5}}&l&rP10{\sqrt{x_7}}&l&u&u\cr}
$$
for $m=2$ and $n=10$ has weight $y_1^2
x_2x_5x_7$. (The weight will always be a product
of  integral powers of the $x_i$ and $y_j$.)

Now I'll define what I call ``dimer tableaux."
Take an $m/2$ by $n/2$ rectangle and split it
into two parts by a path from the lower left
corner to the upper right corner. For example
(with $m=6$ and $n=10$)
$$\Thickness=1pt
\squaresize 15pt
\Young{&&&&r\cr
&&&&lu\cr
d&dr&u&u&\cr}
$$
Then fill in the upper left
part with entries from $1,2,\ldots, n-1$ so that
for adjacent entries 
$\squaresize 10pt\lower 2pt\vbox{\Young{i&j\cr}}$
we have $i<j-1$ and for adjacent entries 
$\squaresize 10pt\lower 8pt\vbox{\Young{i\cr
j\cr}}$ we have $i\le j+1$, and fill in the
lower-right partition with entries from
$1,2\ldots, m-1$ with the reverse inequalities (
$\squaresize 10pt\lower 2pt\vbox{\Young{i&j\cr}}$
implies  $i\le j+1$ and
$\squaresize 10pt\lower 8pt\vbox{\Young{i\cr
j\cr}}$ implies  $i<j-1$). We weight an $i$ in
the upper-left part by $x_i$ and a $j$ in the
lower-right part by $y_j$.

\smallskip\noindent{\bf Theorem:} The sum of the
weights of the $m\times n$ dimer coverings is
equal to the sum of the weights of the $m/2\times
n/2$ dimer tableaux.
\smallskip

My proof is not very enlightening; it essentially
involves showing that both of these are counted
by the same formula.

{\bf Problem 20:} Is there an ``explanation" for
this equality? In particular, is there  a
reasonable bijective proof?  Notes:
\begin{enumerate}
\item[(1)] The case $m=2$ is easy: the $2\times
10$ dimer covering above corresponds to the
$1\times 5$ dimer tableau
$$\squaresize
20pt\Thickness1pt\Young{dx_2&dx_5&drx_7&uy_1&uy_1\cr}$$
(there's only one possibility!)
\item[(2)] If we set $x_i=y_i=0$ when $i$ is even
(so that every two-by-two square of the dimer
covering may be chosen independently), then the
equality is equivalent to the identity
$$\prod_{i,j}(x_i+y_j)=\sum_{\lambda}s_{\lambda}(x)s_{\tilde
\lambda'}(y),$$ (cf.~Macdonald's {\it Symmetric
Functions and Hall Polynomials\/}, p.~37.) This
identity can be proved by a variant of
Schensted's correspondence, so a bijective proof
of the general equality would be essentially a
generalization of Schensted. Several people have
looked at the problem of a Schensted
generalization corresponding to the case in which
$y_i=0$ when
$i$ is even.
\item[(3)] The analogous results in which $m$ or
$n$ is odd are included in the case in which $m$
and $n$ are both even. For example, if we take
$m=4$ and set $y_3=0$, then the fourth row of a
dimer covering must consist of $n/2$ horizontal
dominoes, which contribute
$\sqrt{x_1x_3\cdots x_{n-1}}$ to the weight, so
we are essentially looking at dimer coverings
with three rows.
\end{enumerate}

\end{document}